\documentclass[11pt,reqno]{amsart}
\usepackage{amsfonts,amssymb,amsmath,color}

\DeclareMathOperator{\bmin}{\mathbf{min}}

\newtheorem{theorem}{Theorem}

\newtheorem{lemma}[theorem]{Lemma}
\newtheorem{proposition}[theorem]{Proposition}

\numberwithin{equation}{section} \numberwithin{theorem}{section}

\begin{document}

\title{Pattern Recognition on Oriented Matroids: $\kappa^{\ast}$-Vectors and Reorientations}

\author{Andrey O. Matveev}
\email{andrey.o.matveev@gmail.com}


\begin{abstract}
The components of $\kappa^{\ast}$-vectors associated to a simple oriented matroid $\mathcal{M}$ are the numbers of general or special tope committees for $\mathcal{M}$. Using the principle of inclusion-exclusion, we determine how the reorientations\hfill of\hfill $\mathcal{M}$\hfill on\hfill one-element\hfill subsets\hfill of\hfill its\hfill ground\hfill set\\ affect~$\kappa^{\ast}$-vectors.
\end{abstract}

\maketitle

\pagestyle{myheadings}

\markboth{PATTERN RECOGNITION ON ORIENTED MATROIDS}{A.O.~MATVEEV}

\thispagestyle{empty}

\tableofcontents

\section{Introduction}

Let $\mathcal{M}:=(E_t,\mathcal{T})$ be a simple oriented matroid on the ground set $E_t$ $:=\{1,\ldots,t\}$, with set of topes $\mathcal{T}$; throughout we will suppose that it is {\em simple}, that is, it contains no loops, parallel or {\sl antiparallel\/} elements.

See, e.g.,~\cite{BK,BLSWZ,Bo,DeLoeraRS,RGZ,S,Z} on oriented matroids.

Associated to each element $e\in E_t$ are the corresponding {\em positive half\-space\/}~$\mathcal{T}^+_e:=\{T\in\mathcal{T}:\ T(e)=+\}$ and {\em negative halfspace\/} $\mathcal{T}^-_e:=\{T\in\mathcal{T}:\ T(e)=-\}$ of $\mathcal{M}$.  If $\mathcal{T}_e^{\bullet}\subset\mathcal{T}$ is a halfspace of $\mathcal{M}$ then we denote by $\binom{\mathcal{T}_e^{\bullet}}{j}$ the family of $j$-subsets of the set $\mathcal{T}_e^{\bullet}$.

If $G\subseteq\mathcal{T}$ is a subset of topes then $-G$ stands for the set of their opposites $\{-T:\ T\in G\}$.

If $A\subseteq E_t$ then ${}_{-A}\mathcal{M}$ denotes the oriented matroid obtained from $\mathcal{M}$
by {\em reorientation\/} on the set $A$; if $a\in E_t$ then we write ${}_{-a}\mathcal{M}$ instead of ${}_{-\{a\}}\mathcal{M}$.

A\hfill subset\hfill $\mathcal{K}^{\ast}\subset\mathcal{T}$\hfill is\hfill called\hfill a\hfill {\em tope\hfill committee\/}\hfill for\hfill $\mathcal{M}$\hfill if\hfill for\hfill each\\ element $e\in E_t$ it holds
\begin{equation*}
|\{T\in\mathcal{K}^{\ast}:\ T(e)=+\}|>\tfrac{1}{2}|\mathcal{K}^{\ast}|\ ,
\end{equation*}
see~\cite{M-Halfspaces,M-Layers,M-Existence,M-Three}; in other words, if we replace the components $-$ and $+$ of the maximal covectors of the oriented matroid $\mathcal{M}$ by the real numbers $-1$ and~$1$, respectively, then a collection $\mathcal{K}^{\ast}\subset\mathcal{T}$ is a committee for $\mathcal{M}$ iff the strict inequality
\begin{equation*}
\sum_{T\in\mathcal{K}^{\ast}}T>\pmb{0}
\end{equation*}
holds componentwise.

Let $\mathbf{K}^{\ast}_k(\mathcal{M})$ denote the family of tope committees, of cardinality $k$, for~$\mathcal{M}$, and let $\mathbf{K}^{\ast}(\mathcal{M}):=\dot\bigcup_{1\leq k\leq|\mathcal{T}|-1}\mathbf{K}^{\ast}_k(\mathcal{M})$ denote the family of all tope committees for $\mathcal{M}$. By definition, the $k$th component $\kappa^{\ast}_k(\mathcal{M}):=\#\mathbf{K}^{\ast}_k(\mathcal{M})$ of the vector $\boldsymbol{\kappa}^{\ast}(\mathcal{M})$ $\in\mathbb{N}^{|\mathcal{T}|/2}$, $1\leq k\leq|\mathcal{T}|/2$, is the number of committees in the family $\mathbf{K}^{\ast}_k(\mathcal{M})$.

Similarly, we associate to each family $\overset{\circ}{\mathbf{K}}{}^{\ast}_k(\mathcal{M})$, $1\leq k\leq|\mathcal{T}|/2$, of tope committees, of cardinality $k$, that contain no pairs of opposites, the $k$th component
$\overset{\circ}{\kappa}{}^{\ast}_k(\mathcal{M}):=\#\overset{\circ}{\mathbf{K}}{}^{\ast}_k(\mathcal{M})$
of the vector $\overset{\circ}{\boldsymbol{\kappa}}{}^{\ast}(\mathcal{M})\in\mathbb{N}^{|\mathcal{T}|/2}$.

We always have $\overset{\circ}{\kappa}{}^{\ast}_2(\mathcal{M})=
\kappa^{\ast}_2(\mathcal{M})=0$. The oriented matroid $\mathcal{M}$ is acyclic iff $\overset{\circ}{\kappa}{}^{\ast}_1(\mathcal{M})=
\kappa^{\ast}_1(\mathcal{M})=1$. If $\mathcal{M}$ is not acyclic then $\overset{\circ}{\kappa}{}^{\ast}_1(\mathcal{M})=
\kappa^{\ast}_1(\mathcal{M})=0$ and~$\overset{\circ}{\kappa}{}^{\ast}_3(\mathcal{M})=
\kappa^{\ast}_3(\mathcal{M})$.

If $\mathcal{K}^{\ast}\in\overset{\circ}{\mathbf{K}}{}^{\ast}_{j}(\mathcal{M})$, for some $j$, $1\leq j\leq|\mathcal{T}|/2$, then there are
$|\mathcal{T}|/2 -j$ pairs of topes $\{T,-T\}\subset\mathcal{T}$ such that
$|\mathcal{K}^{\ast}\cap\{T,-T\}|=0$. If we add any such pairs of opposites to the set $\mathcal{K}^{\ast}$ then the resulting set is a committee for $\mathcal{M}$. Thus, given an integer $k$ such that $j\leq k\leq|\mathcal{T}|/2$ and the difference $k-j$ is even, in the family $\mathbf{K}^{\ast}_{k}(\mathcal{M})$ there are exactly $\binom{(|\mathcal{T}|-2j)/2}{(k-j)/2}$ tope committees which contain the committee $\mathcal{K}^{\ast}$ as a subset. We see that
\begin{equation*}
\kappa^{\ast}_k(\mathcal{M})=\sum_{\substack{1\leq j\leq k:\\
j\equiv k\pmod{2}}}\binom{(|\mathcal{T}|-2j)/2}{(k-j)/2}\cdot\overset{\circ}{\kappa}{}^{\ast}_j(\mathcal{M})
\ ,\ \ \ 1\leq k\leq|\mathcal{T}|/2\ ;
\end{equation*}
for example, $\kappa^{\ast}_3(\mathcal{M})=\frac{|\mathcal{T}|-2}{2}\cdot\overset{\circ}{\kappa}{}^{\ast}_1(\mathcal{M})
+\overset{\circ}{\kappa}{}^{\ast}_3(\mathcal{M})$, and
$\kappa^{\ast}_5(\mathcal{M})=\frac{(|\mathcal{T}|-4)(|\mathcal{T}|-2)}{8}$ $\cdot\overset{\circ}{\kappa}{}^{\ast}_1(\mathcal{M})+\frac{|\mathcal{T}|-6}{2}\cdot\overset{\circ}{\kappa}{}^{\ast}_3(\mathcal{M})
+\overset{\circ}{\kappa}{}^{\ast}_5(\mathcal{M})$.

The family $\mathbf{A}^{\ast}(\mathcal{M})$ of {\em anti-committees\/} for the oriented matroid $\mathcal{M}$ is defined as the family $\{-\mathcal{K}^{\ast}:\ \mathcal{K}^{\ast}\in\mathbf{K}^{\ast}(\mathcal{M})\}$.

Let $A$ be any subset of the ground set $E_t$. The tope sets of the oriented matroids ${}_{-A}\mathcal{M}$ and ${}_{-(E_t-A)}\mathcal{M}$ coincide and, thanks to the composite bijection
\begin{gather*}
\mathbf{K}^{\ast}({}_{-A}\mathcal{M})\ \to\ \mathbf{A}^{\ast}({}_{-A}\mathcal{M})\ \to \  \mathbf{A}^{\ast}({}_{-(E_t-A)}\mathcal{M})\ \to \
\mathbf{K}^{\ast}({}_{-(E_t-A)}\mathcal{M})\ ,\\
\mathcal{K}^{\ast}\ \mapsto\ -\mathcal{K}^{\ast} \ \mapsto\ -\mathcal{K}^{\ast}\ \mapsto\ \mathcal{K}^{\ast}\ ,
\end{gather*}
the (anti-)committee structures of ${}_{-A}\mathcal{M}$ and ${}_{-(E_t-A)}\mathcal{M}$ are
identical; in particular, we have
\begin{align*}
\boldsymbol{\kappa}^{\ast}({}_{-A}\mathcal{M})&=\boldsymbol{\kappa}^{\ast}({}_{-(E_t-A)}\mathcal{M})
\intertext{and}
\overset{\circ}{\boldsymbol{\kappa}}{}^{\ast}({}_{-A}\mathcal{M})&=
\overset{\circ}{\boldsymbol{\kappa}}{}^{\ast}({}_{-(E_t-A)}\mathcal{M})\ .
\end{align*}

In this paper we compare $\kappa^{\ast}$-vectors of the oriented matroids $\mathcal{M}$ and~${}_{-A}\mathcal{M}$, where $A:=\{a\}$ are one-element subsets of the ground set~$E_t$.
In~Section~\ref{s:3} we sum up the observations that concern general tope committees and committees containing no pairs of opposites, made in Sections~\ref{s:1} and~\ref{s:2}, respectively.

\section{The Number of Tope Committees}
\label{s:1}

Consider general tope committees for the oriented matroid $\mathcal{M}$ and begin by restating expression~\cite[(3.2)]{M-Halfspaces}:

\begin{lemma}
\label{prop:1}
The number $\#\mathbf{K}^{\ast}_k(\mathcal{M})$ of tope committees, of cardinality~$k$, $1\leq k\leq
|\mathcal{T}|-1$, for the oriented matroid $\mathcal{M}:=(E_t,\mathcal{T})$, is
\begin{equation}
\label{eq:3}
\#\mathbf{K}^{\ast}_k(\mathcal{M})=
\binom{|\mathcal{T}|}{|\mathcal{T}|-\ell} +
\sum_{\substack{\mathcal{G}\subseteq\bigcup_{e\in E_t}\binom{\mathcal{T}^+_e}{\lfloor(\ell+1)/2\rfloor}:\\
1\leq\#\mathcal{G}\leq\binom{\ell}{\lfloor(\ell+1)/2\rfloor},\\
|\bigcup_{G\in\mathcal{G}}G|\leq\ell}}
(-1)^{\#\mathcal{G}}\cdot
\binom{|\mathcal{T}|-|\bigcup_{G\in\mathcal{G}}G|}{|\mathcal{T}|-\ell}\ ,
\end{equation}
where $\ell\in\{k,|\mathcal{T}|-k\}$.
\end{lemma}

Fix an integer $k$, $1\leq k\leq|\mathcal{T}|/2$, a ground element $a\in E_t$, and an integer~$\ell\in\{k,|\mathcal{T}|-k\}$.
If we set
\begin{multline*}
\alpha_k(a,\mathcal{M}):=\binom{|\mathcal{T}|}{|\mathcal{T}|-\ell}+
\sum_{\substack{\mathcal{G}\subseteq\bigcup_{e\in E_t-\{a\}}\binom{\mathcal{T}^+_e(\mathcal{M})}{\lfloor(\ell+1)/2\rfloor}:\\
1\leq\#\mathcal{G}\leq\binom{\ell}{\lfloor(\ell+1)/2\rfloor},\\
|\bigcup_{G\in\mathcal{G}}G|\leq\ell}}
(-1)^{\#\mathcal{G}}\cdot
\binom{|\mathcal{T}|-|\bigcup_{G\in\mathcal{G}}G|}{|\mathcal{T}|-\ell}
\end{multline*}
then, according to~(\ref{eq:3}), we have
\begin{multline}
\label{eq:12}
\kappa^{\ast}_k(\mathcal{M})=\alpha_k(a,\mathcal{M})\\ +
\sum_{\substack{
\mathcal{G}'\;\subseteq\;\binom{\mathcal{T}^+_a(\mathcal{M})}{\lfloor(\ell+1)/2\rfloor}
-\bigcup_{e\in E_t-\{a\}}\binom{\mathcal{T}^+_e(\mathcal{M})}{\lfloor(\ell+1)/2\rfloor}:\ 1\leq\#\mathcal{G}'\leq\binom{\ell}{\lfloor(\ell+1)/2\rfloor},\ |\bigcup_{G\in\mathcal{G}'}G|\leq\ell,\\
\mathcal{G}''\;\subseteq\;\bigcup_{e\in E_t-\{a\}}\binom{\mathcal{T}^+_e(\mathcal{M})}{\lfloor(\ell+1)/2\rfloor}:\ 0\leq\#\mathcal{G}''\leq \binom{\ell}{\lfloor(\ell+1)/2\rfloor}-\#\mathcal{G}',\ |\bigcup_{G\in\mathcal{G}'\dot\cup\mathcal{G}''}G|\leq\ell}}
(-1)^{\#\mathcal{G}'+\#\mathcal{G}''}\\ \cdot
\binom{|\mathcal{T}|-|\bigcup_{G\in\mathcal{G}'\dot\cup\mathcal{G}''}G|}{|\mathcal{T}|-\ell}\ .
\end{multline}
In an analogous expression for $\kappa^{\ast}_k({}_{-a}\mathcal{M})$
the families $\mathcal{G}'$ range over subfamilies  of the family $\binom{\mathcal{T}^-_a(\mathcal{M})}{\lfloor(\ell+1)/2\rfloor}
-\bigcup_{e\in E_t-\{a\}}\binom{\mathcal{T}^+_e(\mathcal{M})}{\lfloor(\ell+1)/2\rfloor}$.

\section{The Number of Tope Committees Containing no Pairs of Opposites}
\label{s:2}

Before proceeding to consider the tope committees that contain no pairs of opposites,
we collect a few observations:

Let $m$ be a positive integer, and $\pm[1,m]$ the $2m$-set
$\{-m,\ldots,-1,1,\ldots,m\}$.
If we fix a subset $W\subseteq\pm[1,m]$ and denote by $-W$ the set $\{-w:\ w\in W\}$ then we have
\begin{multline}
\label{eq:17}
|\!\pm[1,m]|-|W|-2\#\bigl\{\{i,-i\}\subseteq\pm[1,m]:\ |\{i,-i\}\cap W|=0\bigr\}\\=
|W\cup-W|-|W|
\end{multline}
and
\begin{equation}
\label{eq:18}
\#\bigl\{\{i,-i\}\subseteq\pm[1,m]:\ |\{i,-i\}\cap W|=0\bigr\}=
m-\tfrac{1}{2}|W\cup-W|\ .
\end{equation}

Recall that the number of $k$-subsets $V\subset\pm[1,m]$, such that
\begin{equation}
\label{eq:14}
v\in V\ \ \ \Longrightarrow\ \ \  -v\not\in V\ ,
\end{equation}
is $\binom{m}{k}2^k$ --- this is the number of $(k-1)$-dimensional faces of an $m$-dimensional {\em crosspolytope}, see~\cite{HRGZ}.

If $W\neq\pm[1,m]$ then consider some nonempty $k$-set $V\subset\pm[1,m]$ such that~$|V\cap W|$ $=0$ and implication~(\ref{eq:14}) holds. Let $V=V'\dot\cup V''$ be the partition of $V$ into two subsets with the  following properties:
\begin{align}
\label{eq:15}
v'\in V'\ \ \ &\Longrightarrow\ \ \ -v'\in W\ ,\\
\label{eq:16}
v''\in V''\ \ \ &\Longrightarrow\ \ \ -v''\not\in W\ .
\end{align}
Let $|V'|=:j$ and $|V''|=:k-j$, for some $j$. In fact,~(\ref{eq:17}) and~(\ref{eq:18}) imply that\hfill there\hfill are\hfill
$\binom{|W\cup-W|-|W|}{j}$\hfill sets\hfill $V'\subset\pm[1,m]$\hfill such\hfill that\hfill $|V'|=j$,\\ $|V'\cap W|=0$ and~(\ref{eq:15}) holds; there are
$\binom{m-\tfrac{1}{2}|W\cup-W|}{k-j}2^{k-j}$ sets $V''\subset\pm[1,m]$ such that~$|V''|=k-j$, $|V''\cap W|=0$ and~(\ref{eq:16}) holds.

Let $\mathbb{B}(2m)$ denote the Boolean lattice of subsets of the set $\pm[1,m]$. The empty subset of $\pm[1,m]$ is denoted by $\hat{0}$. If $b\in\mathbb{B}(2m)-\{\hat{0}\}$ then
we let $-b$ denote the set of the negations of elements from $b$.

Let $r$ be a rational number, $0\leq r<1$, and $k$ an integer number, $1\leq k$ $\leq m$. If $\varLambda$ is an antichain in $\mathbb{B}(2m)$, such that $\lfloor r\cdot k\rfloor+1\leq\min_{\lambda\in\varLambda}\rho(\lambda)$, then consider the subset
\begin{multline*}
\overset{\circ}{\mathbf{I}}{}_{r,k}\bigl(\mathbb{B}(2m),\varLambda\bigr):=
\bigl\{b\in\mathbb{B}(2m):\\ \rho(b)=k,\ b\wedge-b=\hat{0},\ \rho(b\wedge\lambda)>r\cdot k\ \ \ \forall\lambda\in \varLambda\bigr\}\subset\mathbb{B}(2m)^{(k)}\ ,
\end{multline*}
where\hfill $\rho(\cdot)$\hfill denotes\hfill the\hfill poset\hfill rank\hfill of\hfill an\hfill element\hfill in\hfill $\mathbb{B}(2m)$,\hfill and\hfill $\mathbb{B}(2m)^{(k)}$\\ $:=
\{b\in\mathbb{B}(2m):\ \rho(b)=k\}$. The collection $\overset{\circ}{\mathbf{I}}{}_{r,k}\bigl(\mathbb{B}(2m),\varLambda\bigr)$ is the set of {\em relatively $r$-blocking elements\/} $b\in\mathbb{B}(2m)^{(k)}$ ({\em with the additional property $b\wedge-b=\hat{0}$}) for the antichain $\varLambda$ in the lattice $\mathbb{B}(2m)$; relative blocking is discussed in~\cite{M-Relative}.

Denote by $\mathfrak{I}(\lambda)$ the principal order ideal of the lattice $\mathbb{B}(2m)$ generated by an element $\lambda\in\varLambda$. Using the principle of inclusion-exclusion~\cite{A, St1}, we obtain
\begin{multline}
\label{eq:19}
\bigl|\overset{\circ}{\mathbf{I}}_{r,k}\bigl(\mathbb{B}(2m),\varLambda\bigr)\bigr|=
\binom{m}{k}2^k+
\sum_{D\subseteq\bmin\bigcup_{\lambda\in\varLambda}(\mathbb{B}(2m)^{(\rho(\lambda)-\lfloor r\cdot k\rfloor)}
\cap\mathfrak{I}(\lambda)):\ |D|>0}\\ (-1)^{|D|}\cdot
\sum_{0\leq j\leq k}
\binom{\rho(\bigvee_{d\in D}d\vee-\bigvee_{d\in D}d)-\rho(\bigvee_{d\in D}d)}{j}
\\ \cdot
\binom{m-\frac{1}{2}\rho(\bigvee_{d\in D}d\vee-\bigvee_{d\in D}d)}{k-j}
2^{k-j}\ ,
\end{multline}
where $\bmin\cdot$ denotes the set of minimal elements of a subposet.

Consider the lattice
\begin{equation*}
\mathcal{E}:=\Bigl\{\bigvee_{d\in D}d:\ D\subseteq
\bmin\bigcup_{\lambda\in\varLambda}\bigl(\mathbb{B}(2m)^{(\rho(\lambda)-\lfloor r\cdot k\rfloor)}
\cap\mathfrak{I}(\lambda)\bigr),\ |D|>0 \Bigr\}\ \dot\cup\ \{\hat{0}\}\ ,
\end{equation*}
where $\hat{0}$ is a new least element adjoined. If we let $\mu_{\mathcal{E}}(\cdot,\cdot)$ denote
the {\em M\"{o}bius function\/} of the lattice $\mathcal{E}$, then we have
\begin{multline}
\label{eq:20}
\bigl|\overset{\circ}{\mathbf{I}}_{r,k}\bigl(\mathbb{B}(2m),\varLambda\bigr)\bigr|=
\binom{m}{k}2^k+
\sum_{z\in\mathcal{E}:\ z>\hat{0}}\! \mu_{\mathcal{E}}(\hat{0},z)\\ \cdot
\sum_{0\leq j\leq k}\!
\binom{\rho(z\vee-z)-\rho(z)}{j}\!
\binom{m-\frac{1}{2}\rho(z\vee-z)}{k-j}
2^{k-j}\ ,
\end{multline}
where $\rho(z)$ denotes the poset rank of an element $z$ in the lattice $\mathbb{B}(2m)$.

It was shown in~\cite{M-Halfspaces} that any tope committee $\mathcal{K}^{\ast}\in\mathbf{K}^{\ast}_k(\mathcal{M})$ for the oriented matroid $\mathcal{M}$ is a blocking $k$-set for the family $\bigcup_{e\in E_t}\binom{\mathcal{T}^+_e}{\lfloor(|\mathcal{T}|-k+1)/2\rfloor}$ of tope subsets, of cardinality
$\lfloor(|\mathcal{T}|-k+1)/2\rfloor$, each of which is contained in some positive halfspace, see~Lemma~\ref{prop:1}. As a consequence, the subfamily~$\overset{\circ}{\mathbf{K}}{}^{\ast}_k(\mathcal{M})\subset\mathbf{K}^{\ast}_k(\mathcal{M})$ is precisely the collection of blocking $k$-sets, that are free of opposites, for the family $\bigcup_{e\in E_t}\binom{\mathcal{T}^+_e}{\lfloor(|\mathcal{T}|-k+1)/2\rfloor}$. With the help of~(\ref{eq:19}), we come to the following conclusion:

\begin{lemma}
\label{prop:4}
The number $\#\overset{\circ}{\mathbf{K}}{}^{\ast}_k(\mathcal{M})$ of tope committees, of cardinality~$k$, $1\leq k\leq|\mathcal{T}|/2$, that contain no pairs of opposites, for the oriented matroid~$\mathcal{M}:=(E_t,\mathcal{T})$, is
\begin{multline}
\label{eq:10}
\#\overset{\circ}{\mathbf{K}}{}^{\ast}_k(\mathcal{M})=
\binom{|\mathcal{T}|/2}{k}2^k+
\sum_{\substack{\mathcal{G}\subseteq\bigcup_{e\in E_t}\binom{\mathcal{T}^+_e}{\lfloor(|\mathcal{T}|-k+1)/2\rfloor}:\\
1\leq\#\mathcal{G}\leq\binom{|\mathcal{T}|-k}{\lfloor(|\mathcal{T}|-k+1)/2\rfloor},\\
|\bigcup_{G\in\mathcal{G}}G|\leq|\mathcal{T}|-k}}
(-1)^{\#\mathcal{G}}\\
\cdot
\sum_{0\leq j\leq k}
\binom{|\bigcup_{G\in\mathcal{G}}G\;\cup\;-\bigcup_{G\in\mathcal{G}}G|-|\bigcup_{G\in\mathcal{G}}G|}{j}
\\ \cdot
\binom{\frac{1}{2}(\;|\mathcal{T}|-|\bigcup_{G\in\mathcal{G}}G\;\cup\;
-\bigcup_{G\in\mathcal{G}}G|\;)}{k-j}
2^{k-j}\ .
\end{multline}
\end{lemma}

If $\mathcal{G}$ is a family of tope subsets then we denote by $\boldsymbol{\mathcal{E}}(\mathcal{G})$
the join-semilattice $\{\bigcup_{F\in\mathcal{F}}F:\ \mathcal{F}\subseteq\mathcal{G},\ \#\mathcal{F}>0\}$ that consists of the unions of the sets from the family $\mathcal{G}$ ordered by inclusion and augmented by a new least element $\hat{0}$ which is interpreted as the empty set. The M\"{o}bius function of the lattice~$\boldsymbol{\mathcal{E}}(\mathcal{G})$ is denoted by $\mu_{\boldsymbol{\mathcal{E}}}(\cdot,\cdot)$.

With the help of~(\ref{eq:20}), Lemma~\ref{prop:4} can be restated in the following way:

\begin{proposition}
The number $\#\overset{\circ}{\mathbf{K}}{}^{\ast}_k(\mathcal{M})$ of tope committees which are free of opposites, of cardinality~$k$, $1\leq k\leq|\mathcal{T}|/2$, for the oriented matroid $\mathcal{M}:=(E_t,\mathcal{T})$, is:

\begin{multline*}
\#\overset{\circ}{\mathbf{K}}{}^{\ast}_k(\mathcal{M})
=\binom{|\mathcal{T}|/2}{k}2^k +
\sum_{G\in\boldsymbol{\mathcal{E}}(\bigcup_{e\in E_t}\binom{\mathcal{T}^+_e}{\lfloor(|\mathcal{T}|-k+1)/2\rfloor}):\ 0<|G|\leq|\mathcal{T}|-k}\mu_{\boldsymbol{\mathcal{E}}}(\hat{0},G)\\ \cdot\sum_{0\leq j\leq k}
\binom{|G\;\cup\;-G|-|G|}{j}\!
\binom{\frac{1}{2}(\;|\mathcal{T}|-|G\;\cup\;
-G|\;)}{k-j}2^{k-j}\ .
\end{multline*}
\end{proposition}

If an integer $k$, $1\leq k\leq|\mathcal{T}|/2$, and a ground element $a\in E_t$ are fixed, then we set
\begin{multline*}
\beta_k(a,\mathcal{M}):=
\binom{|\mathcal{T}|/2}{k}2^k+
\sum_{\substack{\mathcal{G}\subseteq\bigcup_{e\in E_t-\{a\}}\binom{\mathcal{T}^+_e(\mathcal{M})}{\lfloor(|\mathcal{T}|-k+1)/2\rfloor}:\\
1\leq\#\mathcal{G}\leq\binom{|\mathcal{T}|-k}{\lfloor(|\mathcal{T}|-k+1)/2\rfloor},\\
|\bigcup_{G\in\mathcal{G}}G|\leq|\mathcal{T}|-k}}
(-1)^{\#\mathcal{G}}\\ \cdot
\sum_{0\leq j\leq k}
\binom{|\bigcup_{G\in\mathcal{G}}G\;\cup\;-\bigcup_{G\in\mathcal{G}}G|-|\bigcup_{G\in\mathcal{G}}G|}{j}
\\ \cdot
\binom{\frac{1}{2}(\;|\mathcal{T}|-|\bigcup_{G\in\mathcal{G}}G\;\cup\;
-\bigcup_{G\in\mathcal{G}}G|\;)}{k-j}
2^{k-j}\ .
\end{multline*}
In view of~(\ref{eq:10}), we have
\begin{multline}
\label{eq:11}
\overset{\circ}{\kappa}{}^{\ast}_k(\mathcal{M})=
\beta_k(a,\mathcal{M})\\+
\sum_{\substack{
\mathcal{G}'\;\subseteq\;\binom{\mathcal{T}^+_a(\mathcal{M})}{\lfloor(|\mathcal{T}|-k+1)/2\rfloor}
-\bigcup_{e\in E_t-\{a\}}\binom{\mathcal{T}^+_e(\mathcal{M})}{\lfloor(|\mathcal{T}|-k+1)/2\rfloor}:\ 1\leq\#\mathcal{G}'\leq\binom{|\mathcal{T}|-k}{\lfloor(|\mathcal{T}|-k+1)/2\rfloor},\ |\bigcup_{G\in\mathcal{G}'}G|\leq|\mathcal{T}|-k,\\
\mathcal{G}''\;\subseteq\;\bigcup_{e\in E_t-\{a\}}\binom{\mathcal{T}^+_e(\mathcal{M})}{\lfloor(|\mathcal{T}|-k+1)/2\rfloor}:\ 0\leq\#\mathcal{G}''\leq \binom{|\mathcal{T}|-k}{\lfloor(|\mathcal{T}|-k+1)/2\rfloor}-\#\mathcal{G}',\ |\bigcup_{G\in\mathcal{G}'\dot\cup\mathcal{G}''}G|\leq|\mathcal{T}|-k}}
(-1)^{\#\mathcal{G}'+\#\mathcal{G}''}\\ \cdot\sum_{0\leq j\leq k}
\binom{|\bigcup_{G\in\mathcal{G}'\dot\cup\mathcal{G}''}G\;
\cup\;-\bigcup_{G\in\mathcal{G}'\dot\cup\mathcal{G}''}G|-|\bigcup_{G\in\mathcal{G}'\dot\cup\mathcal{G}''}G|}{j}
\\ \cdot
\binom{\frac{1}{2}(\;|\mathcal{T}|-|\bigcup_{G\in\mathcal{G}'\dot\cup\mathcal{G}''}G\;\cup\;
-\bigcup_{G\in\mathcal{G}'\dot\cup\mathcal{G}''}G|\;)}{k-j}
2^{k-j}\ .
\end{multline}
In an analogous expression for $\overset{\circ}{\kappa}{}^{\ast}_k({}_{-a}\mathcal{M})$ the families~$\mathcal{G}'$ range over subfamilies of the family $\binom{\mathcal{T}^-_a(\mathcal{M})}{\lfloor(|\mathcal{T}|-k+1)/2\rfloor}$
$-\bigcup_{e\in E_t-\{a\}}\binom{\mathcal{T}^+_e(\mathcal{M})}{\lfloor(|\mathcal{T}|-k+1)/2\rfloor}$.

\section{$\kappa^{\ast}$-Vectors and Reorientations on One-Element Sets}
\label{s:3}

To find the differences of the components of $\kappa^{\ast}$-vectors associated to the oriented matroid $\mathcal{M}$ and to the oriented matroid ${}_{-a}\mathcal{M}$ which is obtained from~$\mathcal{M}$ by reorientation on a one-element subset $\{a\}\subset E_t$, we combine expressions~(\ref{eq:12}) and~(\ref{eq:11}) related to $\mathcal{M}$ with analogous expressions related to~${}_{-a}\mathcal{M}$:

\begin{proposition}
Let $a$ be an element of the ground set $E_t$ of the oriented matroid $\mathcal{M}:=(E_t,\mathcal{T})$. For an integer $k$, $1\leq k\leq|\mathcal{T}|/2$, the sum
\begin{multline*}
\sum_{\substack{
\mathcal{G}''\;\subseteq\;\bigcup_{e\in E_t-\{a\}}\binom{\mathcal{T}^+_e(\mathcal{M})}{\lfloor(|\mathcal{T}|-k+1)/2\rfloor}:\\ 0\leq\#\mathcal{G}''\leq \binom{|\mathcal{T}|-k}{\lfloor(|\mathcal{T}|-k+1)/2\rfloor}-1,\\
|\bigcup_{G\in\mathcal{G}''}G|\leq|\mathcal{T}|-k
}}
(-1)^{\#\mathcal{G}''}\\ \cdot
\Biggl(\ \
\sum_{\substack{
\mathcal{G}'\;\subseteq\;\binom{\mathcal{T}^-_a(\mathcal{M})}{\lfloor(|\mathcal{T}|-k+1)/2\rfloor}
-\bigcup_{e\in E_t-\{a\}}\binom{\mathcal{T}^+_e(\mathcal{M})}{\lfloor(|\mathcal{T}|-k+1)/2\rfloor}:\\ 1\leq\#\mathcal{G}'\leq\binom{|\mathcal{T}|-k}{\lfloor(|\mathcal{T}|-k+1)/2\rfloor}-\#\mathcal{G}'',\\ |\bigcup_{G\in\mathcal{G}'}G|\leq|\mathcal{T}|-k,
}}
(-1)^{\#\mathcal{G}'}\cdot Q(\mathcal{G}',\mathcal{G}'')\\
-
\sum_{\substack{
\mathcal{G}'\;\subseteq\;\binom{\mathcal{T}^+_a(\mathcal{M})}{\lfloor(|\mathcal{T}|-k+1)/2\rfloor}
-\bigcup_{e\in E_t-\{a\}}\binom{\mathcal{T}^+_e(\mathcal{M})}{\lfloor(|\mathcal{T}|-k+1)/2\rfloor}:\\ 1\leq\#\mathcal{G}'\leq\binom{|\mathcal{T}|-k}{\lfloor(|\mathcal{T}|-k+1)/2\rfloor}-\#\mathcal{G}'',\\ |\bigcup_{G\in\mathcal{G}'}G|\leq|\mathcal{T}|-k,
}}
(-1)^{\#\mathcal{G}'}\cdot Q(\mathcal{G}',\mathcal{G}'')\ \;
\Biggr)
\end{multline*}
and the sum
\begin{multline*}
\sum_{\substack{
G''\in\boldsymbol{\mathcal{E}}({\bigcup_{e\in E_t-\{a\}}\binom{\mathcal{T}^+_e(\mathcal{M})}{\lfloor(|\mathcal{T}|-k+1)/2\rfloor})}:\\
0\leq |G''|\leq |\mathcal{T}|-k}}\mu_{\boldsymbol{\mathcal{E}}}(\hat{0},G'')
\\ \cdot
\Biggl(\ \
\sum_{\substack{
G'\in\boldsymbol{\mathcal{E}}(\binom{\mathcal{T}^-_a(\mathcal{M})}{\lfloor(|\mathcal{T}|-k+1)/2\rfloor}
-\bigcup_{e\in E_t-\{a\}}\binom{\mathcal{T}^+_e(\mathcal{M})}{\lfloor(|\mathcal{T}|-k+1)/2\rfloor}):\\ 0<|G'|\leq|\mathcal{T}|-k
}}
\mu_{\boldsymbol{\mathcal{E}}}(\hat{0},G')\cdot \mathfrak{Q}(G',G'')\\
-
\sum_{\substack{
G'\in\boldsymbol{\mathcal{E}}(\binom{\mathcal{T}^+_a(\mathcal{M})}{\lfloor(|\mathcal{T}|-k+1)/2\rfloor}
-\bigcup_{e\in E_t-\{a\}}\binom{\mathcal{T}^+_e(\mathcal{M})}{\lfloor(|\mathcal{T}|-k+1)/2\rfloor}):\\ 0<|G'|\leq|\mathcal{T}|-k
}}
\mu_{\boldsymbol{\mathcal{E}}}(\hat{0},G')\cdot \mathfrak{Q}(G',G'')\ \;
\Biggr)
\end{multline*}
both calculate the difference
\begin{equation*}
\kappa{}^{\ast}_k({}_{-a}\mathcal{M})-\kappa{}^{\ast}_k(\mathcal{M})
\end{equation*}
under
\begin{align*}
Q(\mathcal{G}',\mathcal{G}''):=\binom{|\mathcal{T}|
-|\bigcup_{G\in\mathcal{G}'\dot\cup\mathcal{G}''}G|}{k}\ \ \ \text{and}\ \ \
\mathfrak{Q}(G',G''):=\binom{|\mathcal{T}|
-|G'\cup G''|}{k}\ .
\end{align*}

These sums calculate the difference
\begin{equation*}
\overset{\circ}{\kappa}{}^{\ast}_k({}_{-a}\mathcal{M})-\overset{\circ}{\kappa}{}^{\ast}_k(\mathcal{M})
\end{equation*}
under
\begin{multline*}
Q(\mathcal{G}',\mathcal{G}''):=\sum_{0\leq j\leq k}
\binom{|\bigcup_{G\in\mathcal{G}'\dot\cup\mathcal{G}''}G\;
\cup\;-\bigcup_{G\in\mathcal{G}'\dot\cup\mathcal{G}''}G|-|\bigcup_{G\in\mathcal{G}'\dot\cup\mathcal{G}''}G|}{j}\\
\cdot
\binom{\frac{1}{2}(\;|\mathcal{T}|-|\bigcup_{G\in\mathcal{G}'\dot\cup\mathcal{G}''}G
\;\cup\;-\bigcup_{G\in\mathcal{G}'\dot\cup\mathcal{G}''}G|\;)}{k-j}
2^{k-j}
\end{multline*}
and
\begin{multline*}
\mathfrak{Q}(G',G''):=\sum_{0\leq j\leq k}
\binom{|(G'\cup G'')
\cup-(G'\cup G'')|-|G'\cup G''|}{j}\\
\cdot
\binom{\frac{1}{2}(\;|\mathcal{T}|-|(G'\cup G'')
\cup-(G'\cup G'')|\;)}{k-j}2^{k-j}\ .
\end{multline*}
\end{proposition}


\begin{thebibliography}{10}
\bibitem{A}
M.~Aigner, {\em Combinatorial Theory}, Classics in Mathematics, Reprint of the 1979 original,
Springer-Verlag, Berlin, 1997.

\bibitem{BK}
A.~Bachem and W.~Kern, {\em Linear Programming Duality. An Introduction to Oriented Matroids}. Universitext. Springer-Verlag, Berlin, 1992.

\bibitem{BLSWZ}
A.~Bj\"{o}rner, M.~Las~Vergnas, B.~Sturmfels, N.~White and
G.M.~Ziegler, {\em Oriented Matroids}, Encyclopedia of
Mathematics, {\bf 46}, Cambridge University Press, Cambridge, 1993.
Second edition 1999.

\bibitem{Bo}
J.G.~Bokowski, {\em Computational Oriented Matroids. Equivalence Classes of Matrices within a Natural Framework}. Cambridge University Press, Cambridge, 2006.

\bibitem{DeLoeraRS}
J.A.~De~Loera, J.~Rambau and F.~Santos, {\em Triangulations. Structures for Algorithms and Applications.} Algorithms and Computation in Mathematics, {\bf 25}. Springer-Verlag, Berlin Heidelberg, 2010.

\bibitem{HRGZ}
M.~Henk, J.~Richter-Gebert and G.M.~Ziegler, {\em Basic Properties of Convex Polytopes.} Chapter~16 in {\em   Handbook of Discrete and Computational Geometry} (J.E.~Goodman and J.~O'Rourke, eds.), Chapman~\&~Hall/CRC~Press, Boca Raton, second ed., 2004, 355--382.

\bibitem{M-Halfspaces}
A.O.~Matveev, {\em Pattern Recognition on Oriented Matroids: Halfspaces, Convex Sets and Tope Committees}, {\tt arXiv:1008.4100}.

\bibitem{M-Layers}
A.O.~Matveev, {\em Pattern Recognition on Oriented Matroids: Layers of Tope Committees}, {\tt arXiv:math/0612369}.

\bibitem{M-Existence}
A.O.~Matveev, {\em Pattern Recognition on Oriented Matroids: The
Existence of a Tope Committee}, {\tt arXiv:math.CO/0607570}.

\bibitem{M-Three}
A.O.~Matveev, {\em Pattern Recognition on Oriented Matroids: Three-Tope Committees}, {\tt arXiv:0812.0156}.

\bibitem{M-Relative}
A.O.~Matveev, {\em Relative Blocking in Posets}, J.~Comb.~Optim. {\bf 13} (2007), no~4, 379-403. {\sl Corrigendum}: {\tt arXiv:math.CO/0411026}.

\bibitem{RGZ}
J.~Richter-Gebert and G.M.~Ziegler, {\em Oriented Matroids.} Chapter~6 in {\em   Handbook of Discrete and Computational Geometry} (J.E.~Goodman and J.~O'Rourke, eds.), Chapman~\&~Hall/CRC~Press, Boca Raton, second ed., 2004, 129--152.

\bibitem{S}
F.~Santos, {\em Triangulations of Oriented Matroids}. Mem. Amer. Math. Soc., {\bf 156} (2002), no.~741.

\bibitem{St1}
R.P.~Stanley, {\em Enumerative Combinatorics. Vol. 1}, Corrected reprint of the 1986 original. Cambridge Studies in Advanced Mathematics, {\bf 49}, Cambridge University Press, Cambridge, 1997.

\bibitem{Z}
G.M.~Ziegler, {\em Lectures on Polytopes}. Graduate Texts in Mathematics, {\bf 152}. Springer-Verlag, New~York, 1995.
\end{thebibliography}
\end{document}